%
%
%
%
%
%
%

\scrollmode

%
%
\magnification=\magstep1
\hoffset=1.5cm
\hsize=12cm
\def\pn{\footline={\hss\tenrm\folio\hss}}   

%
%

\def\R{{\rm I\kern-0.2em R\kern0.2em \kern-0.2em}}
\def\N{{\rm I\kern-0.2em N\kern0.2em \kern-0.2em}}
\def\P{{\rm I\kern-0.2em P\kern0.2em \kern-0.2em}}
\def\B{{\rm I\kern-0.2em B\kern0.2em \kern-0.2em}}
\def\C{{\rm C\kern-.4em {\vrule height1.4ex width.08em depth-.04ex}\;}}

\def\ctwo{\C^2}
\def\cn{\C^n}

%
%

\def\ni{\noindent}               
\def\ll{\leftline}
\def\cl{\centerline}

%
%
\def\begin{\ll{}\vskip 5mm \nopagenumbers} 

%
%
\outer\def\beginsection#1\par{\bigskip
  \message{#1}\leftline{\bf\&#1}
  \nobreak\smallskip\vskip-\parskip\noindent}

%
%
\outer\def\proclaim#1:#2\par{\medbreak\vskip-\parskip
    \noindent{\bf#1.\enspace}{\sl#2}
  \ifdim\lastskip<\medskipamount \removelastskip\penalty55\medskip\fi}

%
%
\def\demo#1:{\par\medskip\noindent\it{#1}. \rm}

%
%
\def\endpr{\hfill $\spadesuit$ \medskip}

%
%

%
%
%
%

%
%
%
%

\def\cC{{\cal C}}

\def\cL{{\cal L}}

%
%
%
\def\a{\alpha}

\def\g{\gamma}
\def\d{\delta}
\def\e{\epsilon}
\def\z{\zeta}

\def\l{\lambda}

\def\c{\chi}
\def\t{\theta}

\def\G{\Gamma}

\def\L{\Lambda}

%
%
%
%
\def\bar{\overline}              
\def\bs{\backslash}              
\def\di{\partial}                

%
%
\def\dim{{\rm dim}\,}                    
\def\holo{holomorphic}                   
\def\cont{continuous}                    
\def\nbd{neighborhood}                   
\def\psc{pseudoconvex}                   
\def\ra{real-analytic}                   
\def\psh{plurisubharmonic}               
\def\spsh{strongly\ plurisubharmonic}
\def\ss{\subset\!\subset}                
\def\sd{such that}

\def\iff{if and only if}

\def\phd{proper holomorphic disc}
\def\phm{proper holomorphic map}

\def\dist{{\rm dist}}
\def\te{there exists}
\def\hp{holomorphic polynomial}
\def\st{such that}

\def\eit{e^{i\theta}}

%
%
%
%

\begin
\cl{\bf PROPER HOLOMORPHIC DISCS IN $\C^2$}
\bigskip
\cl{Franc Forstneri\v c and Josip Globevnik}
\bigskip

\rm
\beginsection 1. Results.

Let $U$ denote the open unit disc in $\C$ and $T=bU$ the unit circle.
Let $X$ be a Stein manifold of dimension at least two.
It was proved in [Glo] that for any point $p\in X$ there exists
a \phm\ $f\colon U\to X$ satisfying $f(0)=p$.
We shall call such maps {\it \phd s} in $X$.
For smoothly bounded \psc\ domains in $\C^n$ this was proved
earlier in [FG], and the essential addition in [Glo] was a method
for crossing critical points of a \spsh\ exhaustion function
on $X$. The methods developed in [FG] and [Glo] actually show
the following.

\proclaim 1.1 Theorem: Let $X$ be a Stein manifold
with $\dim X\ge 2$, let $\rho\colon X\to \R$ be
a smooth exhaustion function which is \spsh\ on $\{\rho>M\}$
for some $M\in \R$, and let $d$ be a metric on $X$.
Given a continuous map $h\colon \bar U\to X$ which is holomorphic 
on $U$ and satisfies $\rho(h(\eit))>M$ for $\eit \in T$,
there exists for any pair of numbers $0<r<1$, $\e>0$, and for 
any finite set $A\subset U$ a \phm\ $f\colon U\to X$ satisfying
\item{(i)}   $\lim_{|\z|\to 1 } \rho(f(\z)) = +\infty$,
\item{(ii)}  $\rho(f(\z)) > \rho(h(\z)) -\e$ for $\z \in U$,
\item{(iii)} $d\bigl(f(\z),h(\z)\bigr) < \e$ for $|\z|\le r$,
and
\item{(iv)} $f(\z)=h(\z)$ for  $\z\in A$.

We are interested to what extent does theorem 1.1 hold 
if $\rho$ is a (strongly) \psh\ function 
whose sub-level sets are not necessarily relatively compact.  
Besides its intrinsic interest, we are motivated by the 
question whether it is possible to avoid any closed complex 
hypersurface $L$ in a Stein manifold by \phd s.
Such $L$ is the zero set of a smooth \psh\ function 
$\rho\colon X\to\R_+$ which is \spsh\ on $\{\rho>0\} = X\bs L$;
therefore a positive answer to the first question gives 
\phd s in $X$ avoiding $L$. In this paper we obtain positive 
results in certain model situations in $\C^2$. 
We begin with the following result. 

\proclaim 1.2 Theorem: For each $c<1$ and $M\in \R$ 
the conclusion of theorem 1.1 holds with $X=\C^2$
and the function $\rho_c\colon\C^2\to\R$ given by
$$ \rho_c(z_1,z_2)= \rho_c(x_1+iy_1,x_2+iy_2) =
   x_1^2+x_2^2- c(y_1^2 + y_2^2).     \eqno(1.1)
$$
If on the other hand $c\ge 1$ then for any
\phm\ $f\colon U\to \C^2$ the function $\rho_c\circ f$ 
is unbounded from below on $U$; hence there exist
no proper holomorphic discs satisfying theorem 1.1 (i)
for $\rho=\rho_c$ with $c\ge1$.

Note that $\rho_c$ is \spsh\ if $c<1$, strongly plurisuperharmonic
if $c>1$, and $\rho_1(z_1,z_2)=\Re(z_1^2+z_2^2)$ is pluriharmonic.

The second statement in theorem 1.2 (for $c\ge 1$) 
can be seen by applying theorem 1.5 (d) below 
to the function $g=f_1^2+f_2^2$: since its range at any 
boundary point $\eit\in T$ omits at most a polar set
in $\C$, its real part $\Re g=\rho_1(f_1,f_2)$ is unbounded 
from below. Since $\rho_c\le \rho_1$ for $c\ge 1$, the same is 
true for $\rho_c\circ f$. The first part of theorem 1.2 
(for $c<1$) is proved in section 3. 

When $c>0$, $\rho_c$ is not an exhaustion function on $\C^2$. 
For $0<c<1$ theorem 1.2 gives \phm s $f\colon U\to\C^2$ 
with images $f(U)$ contained in the real cone $\G_c = \{\rho_c>0\}$
with axis $\R^2=\{y=0\}$. Moreover, when $c>1$ 
we can apply theorem 1.2 with 
$- \rho_c(z)/c = y_1^2+y_2^2 - {1\over c} (x_1^2 +x_2^2)$
to obtain a \phm\ $f\colon U\to \C^2$ whose image 
avoids $\Gamma_c$. {\it This gives proper holomorphic discs
in $\C^2$ avoiding relatively large real cones.}
On the other hand, no proper holomorphic disc (in fact, 
no transcendental complex curve) in $\C^2$ can avoid a 
nonempty open complex cone; see theorem 2 in [SW] 
and theorem 1.5 below.

Our next result concerns discs avoiding pairs of complex lines
in $\C^2$.

\proclaim 1.3 Theorem: There exists a \phm\
$f=(f_1,f_2)\colon  U\to\ctwo$ whose
image $f(U)$ is contained in $(\c*)^2=\C^2\bs\{zw=0\}$.

Writing $f\colon  U\to (\c*)^2$ as
$f=\bigl(e^{g_1}, e^{g_2} \bigr)=
   \bigl(e^{u_1+iv_1},e^{u_2+iv_2}\bigr),
$
we have $|f|^2=|f_1|^2+|f_2|^2 = e^{2u_1}+e^{2u_2}$, and $f$ is
proper as a map into $\ctwo$ \iff\ $\max\{u_1,u_2\}$ tends to
$+\infty$ at the boundary of $U$.
Thus theorem 1.3 is equivalent to

\pn

\proclaim 1.4 Theorem: There exists a pair of harmonic
functions $u_1,u_2$ on the disc $U$ such that
$$
   \lim_{|\z|\to1} \max\{u_1(\z),u_2(\z)\} = +\infty.
$$

Theorem 1.3 is a special case of theorem 4.1 in section 4
below. A different proof of theorem 1.4 was shown to us
by J.-P.\ Rosay (private communication).

It would be interesting to know whether proper discs
in $\C^2$ can avoid any given finite collection of complex lines.
Part (d) in theorem 1.5 shows that such a disc cannot avoid a
non-polar set of complex lines through the origin
(or parallel complex lines) in $\C^2$. The same holds if 
we replace the disc by any transcendental complex curve 
(Sibony and Wong [SW], Theorem 2).
H.\ Alexander [Ale] proved in 1975 that for parallel
lines in $\C^2$ this is the only obstruction:
\it  If $E \subset \C$ is a closed polar set containing
at least two points, there exists a \phm\
$f=(f_1,f_2)\colon U\to \C^2$ such that
$f_1 \colon U\to \C\bs E$ is a universal
covering map of the disc onto $\C\bs E$. \rm
We don't know whether an analogue of Alexander's
result holds for complex lines through the origin.


In the remainder of this section we discuss
the boundary behavior of proper holomorphic maps
$f=(f_1,f_2)\colon U\to\C^2$ at the circle $T=\{|\z|=1\}$.
We must recall some basic
notions from the theory of cluster sets of meromorphic
functions on the disc; we refer to Chapter 8
in the monograph [CL] (see section 5 below
for more details).

Let $g$ be a meromorphic function on $U$. A point $\eit\in T$ at
which the (unrestricted) cluster set of $g$ equals
$\bar {\C}=\C\cup\{\infty\}$ is called a {\bf Weierstrass point}
of $g$. If the restricted cluster set of $g$ at $\eit$
within each conical region in $U$ with vertex $\eit$ equals
$\bar \C$ then $\eit$ is called a {\bf Plessner point} of $g$.
A point $\eit$ at which $g$ has a non-tangential limit
(a limit as $\z\to\eit$ within any cone in $U$
with vertex $\eit$) is called a {\bf Fatou point} of $g$,
and the set of all Fatou point is the {\bf Fatou set} of $g$.
The {\bf range} of $g$ at $\eit$, denoted
$R(g,\eit)$, consists of all $\a\in \bar \C$
such that $g(\z_j)=\a$ for points in a sequence 
$\z_j\in U$ with $\lim_{j\to\infty} \z_j=\eit$.

%
%
\proclaim 1.5 Theorem:
Let $f=(f_1,f_2)\colon  U\to\ctwo$ be a proper holomorphic map
of the disc to $\C^2$. Let $P,Q$ be nonconstant \holo\
polynomials on $\ctwo$ whose leading order homogeneous
parts have no common divisor. Denote by $g$ any of the
following (meromorphic) functions:
(i) $f_1$ or $f_2$, (ii) $f_1/f_2$, (iii) $P(f_1,f_2)$,
(iv) $P(f_1,f_2)/Q(f_1,f_2)$. Then
\item{(a)} the Fatou set of $g$ has Lebesgue measure zero in $T$,
\item{(b)} every point of $T$ is a Weierstrass point of $g$,
\item{(c)} almost every point of $T$ is a Plessner point of $g$, and
\item{(d)} for every $\eit\in T$ the set
$\C\bs R(g,\eit)$ is polar.

Theorem 1.5 is proved in section 5. Part (d) can be interpreted
as a result on polynomial hulls as follows. We define the
polynomial hull $\widehat{K}$ of an arbitrary
subset $K\subset \cn$ as the intersection of all closed
set in $\cn$ of the form $\{\Re P\le 0\}$ containing $K$,
where $P$ is a holomorphic polynomial. For compact sets
this coincides with the usual definition of the polynomial hull.
Clearly $\widehat K$ is contained in the closed convex hull
of $K$. Theorem 1.5 (d) implies

\proclaim 1.6 Corollary:
If $f\colon  U\to\ctwo$ is a \phm\ then for
each open set $D\subset \C$ intersecting $T$
the polynomial hull of $f(U\cap  D)$ equals $\ctwo$
(and hence its closed convex hull also equals $\ctwo$).

Theorem 1.5 does not generalize directly to proper
maps $f\colon U\to\C^n$ for $n>2$. Namely, if
$(f_1,f_2)\colon U\to\C^2$ is proper holomorphic
and if $f_3$ is any holomorphic function on $U$ then
$(f_1,f_2,f_3)\colon U\to\C^3$ is also proper holomorphic;
thus the addition of the third component need not
enlarge the cluster set at any boundary point.

In the Appendix we comment on the proof of theorem
1.1 in [Glo]. Let $\rho\colon X\to\R$ be  a \spsh\ Morse 
exhaustion function on a Stein manifold $X$ of dimension $\ge2$. 
We show that one can push the boundary of an analytic disc 
in $X$ over a critical level of $\rho$ by using the gradient 
flow of $\rho$. This creates a non-holomorphic contribution
which can be cancelled off during a later stage
of the lifting procedure (this was the crucial
observation in [Glo]).

\demo Acknowledgements:
We thank D.\ Marshall who showed us a simple proof of 
Frostman's theorem to the effect that a meromorphic
function on $U$ which omits a set of positive capacity
has bounded Nevanlinna characteristic (and hence its 
Fatou set has full measure in $T$; see section 5). 
We also thank N.\ Sibony for pointing out the reference [SW].
Research of the first author was supported in part by an NSF grant,
by the Vilas foundation at the University of Wisconsin--Madison,
and by the Ministry of Science of the Republic of Slovenia.
Research of the second author was supported in part by
the Ministry of Science of the Republic of Slovenia.

%
%
%
%
\beginsection 2. Lifting holomorphic discs.

In this section we describe a general method for
lifting the boundary of an analytic disc in $\cn$
to a higher level set of a \spsh\ function $\rho\colon\cn\to\R$.
This method was developed in [FG], but for our present needs
we need more precise estimates for the amount of
possible lifting at each step of the process.

\proclaim 2.1 Proposition: Let $\l\colon T\times \bar U \to\cn$
be a continuous map such that for each $\z\in T$ the map
$\l_\z =\l(\z,\cdotp) \colon \bar U\to \cn$ is holomorphic in $ U$
and $\l_\z(0)=0$. Given numbers $\e>0$ and $0<r<1$, there exists
a holomorphic polynomial map $h\colon \C \to\cn$ satisfying
\item{(i)}  $\dist\bigl( h(\z),\l_\z(T)\bigr) < \e$ \ \ $(\z\in T)$,
\item{(ii)} $\dist\bigl( h(t\z),\l_\z(\bar U)\bigr) < \e$ \ \
$(\z\in T,\ r\le t\le 1)$, and
\item{(iii)} $|h(\z)| < \e$\ \ $(|\z|\le r)$.

\demo Proof: It suffices to show that $\l$ can be approximated uniformly
on $T\times\bar U$ by maps of the form
$$ \tilde\l(\z,w) =  {w\over \z^M} \sum_{j=1}^N A_j(\z) w^{j-1}, \eqno(2.1) $$
where the $A_j$'s are holomorphic polynomials and $M,N$ are positive integers.
The polynomial map
$$
  h(\z)=\tilde\l(\z,\z^K) = \z^{K-M} \sum_{j=1}^N A_j(\z) \z^{(j-1)K}
$$
then satisfies proposition 2.1 provided that the approximation of $\l$
by $\tilde \l$ is sufficiently close and the integer $K\ge M$ is chosen
sufficiently large.

We begin by replacing $\l$ by $(\z,w)\mapsto \l(\z,sw)$ for a suitable
$s<1$ sufficiently close to $1$. Denoting the new map again by $\l$
we may thus assume that $\l_\z$ is holomorphic in a larger disc $|w|<1/s$
for each $\z\in T$. We expand $\l$ in Taylor series with respect to $w$
and approximate it uniformly on $b U\times T$ by a Taylor
polynomial $\l_N(\z,w)= \sum_{j=1}^N a_j(\z) w^j$ with continuous
coefficients $a_j \colon T\to \cn$. (The coefficient $a_0$ is
zero since $\l(\z,0)=0$.) Finally we approximate each $a_j$ uniformly
on $T$ by a map $A_j(\z)/\z^M$ for some holomorphic polynomial
$A_j$ and some integer $N$ which can be chosen to be independent of $j$.
This gives the desired approximation of $\l$ by a map of the form
(2.1).
\endpr

\proclaim 2.2 Corollary: Let $g_0\colon \bar U \to\cn$ be a continuous
map that is holomorphic in $ U$ and let $\l$ be as in proposition 2.1.
Suppose that $\rho\colon \cn \to\R$ is a real continuous function
such that for some constants $C_0<C_1$ and $0<r<1$ we have
\item{(a)} $\rho\bigl( g_0(\z) + \l(\z,w) \bigr) = C_1$\ \
   $(\z \in T,\ w\in T)$,
\item{(b)} $\rho\bigl( g_0(\z) + \l(\z,w) \bigr) > C_0$\ \
   $(\z \in T,\ w\in \bar U)$, and
\item{(c)} $\rho\bigl(g_0(\z)\bigr) > C_0$\ \ $(r\le |\z| \le 1)$.

\ni\sl Then for each $\e>0$ there exists a holomorphic polynomial map
$g\colon \C\to\cn$ satisfying
\item{(i)} $\big| \rho(g(\z)) - C_1 \big| <\e$\ \ $(\z\in T)$,
\item{(ii)} $\rho (g(\z)) > C_0$\ \ $(r\le |\z| \le 1)$, and
\item{(iii} $\big|g(\z)-g_0(\z)\big| < \e$\ \ $(|\z|\le r$).

\demo Proof: Take $g(\z)=\tilde g_0(\z)+ h(\z)$, where $\tilde g_0$ is
a polynomial approximation of $g_0$ and $h$ is a suitably chosen map provided
by proposition 2.1.
\endpr

Assume now that $\rho\colon \cn \to \R$ is a function of  class $\cC^2$.
For each fixed $z$ we write
$$
   \rho_z(w)= \rho(z+w)-\rho(z) =  \Re Q_z(w) + \cL_z(w)+o(|w|^2), \eqno(2.2)
$$
where
$$ \eqalign{
   Q_z(w) &= 2\sum_{j=1}^n {\di\rho\over\di z_j}(z)w_j
           + \sum_{j,k=1}^n {\di^2 \rho \over \di z_j \di z_k}(z) w_j w_k \cr
   \cL_z(w) &= \sum_{j,k=1}^n {\di^2 \rho \over \di z_j \di \bar z_k}(z)
                   w_j \bar w_k \cr}
$$
($Q_z$ is the \it Levi polynomial \rm and $\cL_z$ is the \it Levi form \rm
of $\rho$ at $z$). The set
$$
   \Lambda_z = \{w\in \cn\colon Q_z(w)=0 \}   		 \eqno(2.3)
$$
is a quadratic complex hypersurface in $\cn$ and we have
$\rho_z(w)= \cL_\rho(z;w) + o(|w|^2)$ for $w\in \Lambda_z$.  
For $c>0$ we denote by $B(z;c)$ the connected component of
the sublevel set $\{w\in \L_z\colon \rho_z(w) < c \}$
which contains the point $0\in \L_z$. If $\rho$ is \spsh\ near $z$
(i.e., its Levi form $\cL_z$ at $z$ is positive definite) and
if $\di\rho(z)\ne 0$ (so that the hypersurface $\L_z$ is smooth
near $0$), then for all sufficiently small $c>0$  the set
$B(z;c)$ is diffeomorphic to the real $(2n-2)$-dimensional ball.
Moreover, if $C>0$ is such that the function $\rho_z|\L_z$ has
no critical points on $B(z;C)$ other than the point $0$,
Morse theory shows that for $0< c \le C$ the sets $B(z;c)$
are complex manifolds diffeomorphic to the $(2n-2)$-ball.
(We include the singularities of $\L_z$ among the critical
points of $\rho_z|\L_z$.) In particular, when $n=2$, these sets are
complex one-dimensional and hence conformally equivalent to the disc.
We state the next proposition only for $n=2$ since we shall only
need this case.

\proclaim 2.3 Proposition: Let $g_0\colon \bar U \to\ctwo$ be a continuous
map that is holomorphic in $U$ and let $\rho \colon \C^2\to \R$
be a $\cC^2$ function which is \spsh\ in a \nbd\ of $g_0(T)$
and has no critical points on $g_0(T)$.
Suppose that $C\colon T\to (0,\infty)$ is a continuous
function \st\ the function $\rho_{g_0(\z)}|\L_{g_0(\z)}$ (2.2) has no
critical points on $B(g_0(\z);C(\z)) \bs\{0\}$ for each $\z\in  T$.
Then for each $\e>0$ and $0<r<1$ there is polynomial map
$g\colon \C\to\ctwo$ satisfying
\item{(i)}    $|\rho(g(\z)) - \rho(g_0(\z))- C(\z)| < \e$\ \ $(\z\in T)$,
\item{(ii)}   $\rho(g(\z)) > \rho(g_0(\z))-\e$\ \  $(\z\in\bar U)$, and
\item{(iii)}  $|g(\z)-g_0(\z)| < \e$\ \ $(|\z|\le r)$.

\demo Proof:
We have seen above that for each $\z\in T$ the set
$B(g_0(\z);C(\z)) \subset \L_{g_0(\z)}$ is conformally equivalent
to the disc $U$. Decreasing $C(\z)$ slightly (so that $B(g_0(\z);r)$
is still biholomorphic to $U$ for some $r>C(\z)$) we can
obtain a parametrization
$\l_\z\colon \bar U\to \bar B(g_0(\z);C(\z))$ ($\z\in T$),
depending continuously on $(\z,w)\in T\times \bar U$,
such that $\l_\z$ is \holo\ in $U$ and $\l_\z(0)=g_0(\z)$
for each $\z\in T$. The result now follows from
proposition 2.1 applied to the family of discs $\l_\z$.
\endpr

If $K_0 \ss K_1\ss \ctwo$ is a pair of compact sets
such that $\rho$ is \spsh\ and has no critical points on
$K_1$, there is a constant $C>0$ such that $\rho_z|\L_z$ has no
critical points on $B(z;C)\bs\{0\}$ for each $z\in K_0$.
Hence proposition 2.3 provides a uniform lifting of the
boundary of an analytic disc (with respect to $\rho$) as 
long as the boundary remains in $K_0$. If the set
$A(c_0,c_1)= \{x\in X\colon c_0\le \rho(x)\le c_1\}$ is 
compact for some $c_0<c_1$ and if $\rho$ is \spsh\ and
without critical points on this set, proposition 2.3 allows
us to lift the boundary of an analytic disc in $X$ from the 
level $\rho=c_0$ to the level $\rho=c_1$. Unfortunately
this breaks down in general if the level sets of $\rho$ 
are not compact. In this case we need a more precise analysis 
which we shall do for the function (1.1).

\proclaim 2.4 Proposition: Let $\rho_c$ be the function (1.1).
If $c<1$ there exists a number $a=a(c)>0$ with the following property:
For each continuous map $h\colon \bar U\to\ctwo$, holomorphic in $U$,
such that 
$ m(h) = \inf\{\rho_c(h(\z)) \colon |\z|=1\} >0$,
and for each pair of numbers $\e>0$ and $0<r<1$ there exists a \holo\
polynomial map $g\colon \C\to\ctwo$ satisfying
\item{(i)}   $m(g)\ge (1+a) m(h)$,
\item{(ii)}  $\rho_c(g(\z)) > \rho_c(h(\z)) - \e$\ \ $(|\z|\le 1)$, and
\item{(iii)} $|g(\z)-h(\z)|<\e$\ \ $(|\z|\le r)$.

\demo Proof:  Note $\rho_c$ is \spsh\ when $c<1$.
Fix such a $c$ and write $\rho=\rho_c$. The only critical point
of $\rho$ is $z_1=z_2=0$. Proposition 2.4 follows
immediately from proposition 2.3 and the following

\proclaim 2.5 Lemma:  Let $\rho=\rho_c$ for some $c<1$
be given by (1.1). There is a constant $a=a(c)>0$ such that for
each $z\in\ctwo$ with $\rho(z)>0$ the function $\rho_z|\L_z$
has no critical points on $B(z;a\rho(z)) \bs \{0\}$.
\rm

\demo Proof: A calculation shows that $\rho_z(w)= \Re Q_z(w) + \cL_z(w)$, where
$$ \eqalign{
    Q_z(w)   &=  2(x_1+icy_1)w_1+2(x_2+icy_2)w_2
                 + {1\over 2}(1+c)\bigl( w_1^2 + w_2^2\bigr) \cr
    \cL_z(w) &= {1\over 2}(1-c)(|w_1|^2+|w_2|^2) = {1\over 2}(1-c)|w|^2. \cr}
$$
It suffices to consider the case $0<c<1$. If 
$w\in \Lambda_z$ then
$$
    \rho_z(w)=\rho(z+w)-\rho(z)= {1\over 2}(1-c)|w|^2    \eqno(2.4)
$$
The critical points of $\rho_z|\L_z$ are precisely those points 
$w\in \L_z$ at which the complex gradients $\di Q_z$ 
and $\di \rho_z$ (with respect to the variable $w=(w_1,w_2)\in\ctwo$)
are $\C$-linearly dependent. 
This set will include any singular points of $\L_z$. 
By (2.4) we may replace $\di\rho_z$ by $\di |w|^2$. Set
$h(x+iy)= x+icy$, so $|h(x+iy)|^2= x^2+c^2y^2$.
We have
$$ \di Q_z(w) = \bigl( 2h(z_1)+(1+c)w_1,2h(z_2)+(1+c)w_2 \bigr),
   \qquad \di |w|^2 =  (\bar w_1,\bar w_2).
$$
This gives the following system of two equations for $w$, in which
the first is the colinearity equation between $\di Q_z$ and 
$\di |w|^2$ (after conjugation) and the second is 
$Q_z(w)=0$:
$$ \eqalignno{
   2\bar{h(z_2)} w_1 - 
   2\bar{h(z_1)} w_2 &= -(1+c)(w_1\bar w_2 - \bar w_1 w_2) \cr
   4h(z_1) w_1 + 4h(z_2)w_2 &= -(1+c)(w_1^2+w_2^2). &(2.5) \cr}
$$
It suffices to obtain a good lower estimate for the norm $|w|$ of 
any nonzero solution of (2.5) in terms of $|z|$. 
We apply Cramer's formula to express $w_1$ and $w_2$
from the linear part in terms of the right hand side terms in (2.5).
The determinant of the matrix of coefficients is
$W(z)= 8 \bigl( |h(z_1)|^2 + |h(z_2)|^2\bigr) \ge c'|z|^2$
where $c'>0$ depends only on $c$. If we replace one of the 
columns of the coefficient matrix by the right hand side 
then each term in the corresponding determinant is of the form
constant times $h(z_j)w_kw_l$ for some $j,k,l\in\{1,2\}$.
Hence we can estimate the determinant from above by the Cauchy-Schwarz
inequality and thus obtain the following estimate for the solutions of (2.5):
$$ |w_j| \le {c_2 \bigl( |h(z_1)|^2 + |h(z_2)|^2\bigr)^{1/2} 
   |w|^2 \over W(z)} \le {c_3 |w|^2 \over |z|} \quad (j=1,2). 
$$
This gives $|w|\le c_4 |w|^2/|z|$ and therefore $|w|\ge c_5|z|$ for
any nonzero solution $w$ of (2.5), where $c_5>0$ depends
only on $c$. Since $w\in\L_z$, (2.4) gives
$$ \rho(z+w)\ge \rho(z)+c_6|z|^2 \ge \rho(z)+c_7\rho(z) $$
for some $c_7>0$. Thus any constant $a<c_7$ satisfies 
lemma 2.5.
\endpr

%
%
%
%
\beginsection 3. Proper discs in cones in $\ctwo$ with real axis.

In this section we prove theorem 1.2. If the constant $M$ in 
the theorem is negative, we first apply the procedure described
in [Glo] to cross the critical point of $\rho_c$ at $(0,0)$ 
and thus push the boundary of the given initial analytic
disc $h$ to the set $\rho_c>0$ while changing $h$ as little 
as desired on $\{|\z|\le r\}$. Hence it suffices to prove 
theorem 1.2 for $M\ge 0$. In this case the result follows 
immediately from the following.

\proclaim 3.1 Theorem: Let $c<1$, $M\ge 0$, and let
$\rho=\rho_c$ be the function (1.1). Given a \cont\ map
$h\colon  \bar U\to \ctwo$, \holo\ in $U$, such that
$\rho(h(\z)) > M$ for $|\z|=1$, there exists for each
$\e>0$ and $0<r_1<1$ a \phm\ $f\colon  U\to\ctwo$ satisfying
\item{(i)}   $\lim_{|\z|\to 1} \rho_c(f(\z)) = +\infty$,
\item{(ii)}  $\rho(f(\z)) > \rho(h(\z)) -\e$\ \ $(|\z|<1)$, and
\item{(iii)} $|f(\z)-h(\z)| < \e$\ \ $(|\z|\le r_1)$.

\demo Proof:
It suffices to consider the case $0<c<1$. Fix numbers
$M>0$, $0<r<1$, $\e>0$ and a map $h$ as in the statement of
theorem 3.1 and write $M_1=M$, $\e_1=\e$, $f_1=h$. Let
$a>0$ be the number given by proposition 2.4 for the pair
$c$ and $M_1$. Set
$$
    M_k=(1+a)^{k-1}M_1,\quad \e_k=\e/2^{k-1},\quad   k=2,3,4,\ldots
$$
We inductively construct a sequence of polynomial maps
$f_k\colon \bar U \to\ctwo$ and a sequence of numbers
$0<r_1<r_2<r_3<\ldots<1$ with $\lim_{k\to\infty} r_k=1$
such that the following hold for each $k\ge2$:
\item{(a$_k$)} $\rho(f_k(\z)) > M_k$\ \ $(r_k\le |\z|\le 1)$,
\item{(b$_k$)} $\rho(f_k(\z)) > \rho(f_{k-1}(\z)) - \e_{k-1}$
\ \ $|\z|\le 1)$, and
\item{(c$_k$)} $|f_k(\z)-f_{k-1}(\z)| <\e_{k-1}$\ \ $(|\z|\le r_{k-1})$.

\medskip
The construction proceeds as follows. By assumptions the
condition (a$_1$) holds for $|\z|=1$. By continuity we can
increase $r_1$ such that (a$_1$) holds for $r_1\le |\z|\le 1$.
Proposition 2.4 gives a map $f_2$ such that $\rho(f_2(\z)) >M_2$
for $|\z|=1$ and such that (b$_2$) and (c$_2$) hold. By continuity
we can choose a number $r_2<1$ sufficiently close to $1$
such that (a$_2$) holds for $r_2\le |\z|\le 1$.

This process can be continued inductively. If we already have
$f_{k-1}$, proposition 2.4 gives the next map $f_k$ which satisfies
(a$_k$) initially only for $|\z|=1$, and it satisfies (b$_k$) and
(c$_k$). By continuity we can choose $r_k<1$ sufficiently close to
$1$ so that (a$_k$) holds. We can thus 
insure that $\lim_{k\to\infty} r_k =1$.

Condition (c) insures that $f=\lim_{k\to\infty} f_k \colon U\to \ctwo$
exists uniformly on compacts in $U$. For $|\z|\le r_1$ we have
$$ |f(\z)-f_1(\z)| \le
   \sum_{k=1}^\infty |f_{k+1}(\z)-f_k(\z)| < \sum_{k=1}^\infty \e_{k+1} =\e. $$
This proves (iii) since $h=f_1$. For a fixed $\z\in U$ and $k\ge1$ we have
$$ \eqalign{ \rho(f(\z)) &= \lim_{j\to\infty} \rho(f_j(\z)) 
    = \rho(f_k(\z)) + \sum_{j=k}^\infty \bigl(
     \rho(f_{j+1}(\z)) - \rho(f_j(\z)) \bigr)  \cr
    &> \rho(f_k(\z)) - \sum_{j=k}^\infty \e_{j+1} 
	= \rho(f_k(\z)) - \e_k. \cr}
$$
For $k=1$ we get (ii) in the theorem. For points $\z$ in the annulus
$r_k\le |\z|<1$ we get
$\rho(f(\z)) > \rho(f_k(\z)) - \e_k > M_k-\e$.
Since $\lim_{k\to\infty} M_k =\infty$, this implies (i) 
and completes the proof of theorem 3.1.
\endpr

%
%
%
%
\beginsection 4. Proper discs in $\ctwo$ which omit a pair of lines.

Theorem 1.3 follows from the following more precise result.

\proclaim 4.1 Theorem:
Let $n\ge2$. Given a \cont\ map $h=(h_1,h_2,\ldots,h_n)\colon  \bar U\to \cn$
which is \holo\ in $U$ and given a number $0<r<1$ \sd\ the components $h_j$
have no zeros in $\{\z\colon r\le |\z|\le 1\}$, there exists for each $\e>0$
a \phm\ $f=(f_1,f_2,\ldots,f_n)\colon  U\to\cn$ \sd\ the $f_j$'s  have no
zeros in $\{\z\colon r\le |\z|<1\}$ and $|f(\z)-h(\z)|<\e$ for $|\z|\le r$.

We shall give details only for $n=2$. By factoring out the 
(finitely many) zeros of the $h_j$'s we can reduce to the case when 
the $h_j$'s have no zeros on $\bar U$. We seek a solution in the form
$f=\bigl(e^{g_1}, e^{g_2} \bigr)= \bigl(e^{u_1+iv_1},e^{u_2+iv_2}\bigr)$                                                      
for some holomorphic map $g=(g_1,g_2) \colon  U\to\ctwo$. Set
$$
   \rho(x_1+iy_1,x_2+iy_2) = \max\{x_1,x_2\}.                 \eqno(4.1)
$$
Since $|f|^2=|f_1|^2+|f_2|^2 = e^{2u_1}+e^{2u_2}$, $f$ is proper
into $\ctwo$ \iff\ $\rho(g(\z)) = \max\{u_1(\z),u_2(\z)\}$ tends
to $+\infty$ as $|\z|\to 1$. Such map $g$ will be obtained as the
limit $g=\lim_{k\to\infty} g_k$ of an inductively constructed sequence
$g_k$, where the inductive step from $g_{k-1}$ to $g_k$ will be
furnished by corollary 2.2. To this end we need a suitable family 
of analytic discs which we now construct.

\proclaim 4.2 Proposition:  Let $\rho$ be the function (4.1). Given a
compact set $K\ss \ctwo$ and constants $C_0,C_1\in \R$ such that
$C_0 < \rho(z) < C_1$ ($z\in K$), there is a continuous map
$\l\colon K\times  \bar U \to\ctwo$ such that for each $z\in K$
the map $\l(z,\cdotp)\colon  U \to \ctwo$ is holomorphic and 
\item{(i)}   $\rho(\l(z,w)) =C_1$\ \ $(z\in K,\ |w|=1)$, 
\item{(ii)}  $\rho(\l(z,w)) > C_0$\ \ $(z\in K,\ |w|\le 1)$.

\demo Proof: We follow the proof of Bochner's tube theorem 
(see [H\"or], p.\ 41). We first describe the model situation. 
Write the coordinates on $\ctwo$ in the form $z=x+iy$, with
$x,y\in \R^2$, and identify $\R^2$ with $\{y=0\}\subset \ctwo$. Set
$$ \eqalign{ k &= \{(x_1,0)\colon 0\le x_1\le 1\} \cup
                 \{(0,x_2)\colon 0\le x_2\le 1\}  \cr
          K_\e &= \{x+iy\in\ctwo\colon x\in k,\ |y|^2\le 1/\e\} \cr
         co(k) &= \{(x_1,x_2)\colon x_1\ge0,\ x_2\ge0,\ x_1+x_2\le 1\} \cr
         \g_\e &= \{(x_1,x_2)\in co(k)\colon x_1+x_2-\e(x_1^2+x_2^2)= 1-\e \}  \cr
         \G_\e &= \{(z_1,z_2)\in \ctwo \colon (x_1,x_2)\in co(k),\
                    z_1+z_2-\e(z_1^2+z_2^2)= 1-\e \}  \cr}
$$

\proclaim 4.3 Lemma: (Notation as above)\ \ There is an $\e_0>0$ such that
for each $\e$ with $0<\e<\e_0$ the set $\G_\e$ is a holomorphic disc
with boundary contained in $K_\e$, $\G_\e\cap \R^2=\g_\e$, and 
$\g_\e$ is a smooth real-analytic curve contained in the convex hull 
$co(k)$ of $k$. The union $\bigcup_{0<\e<\e_0} \g_\e$ contains every point 
in the interior of $co(k)$ and sufficiently close to the open segment
$\g_0= \{(x_1,1-x_1)\colon 0< x_1 <1\}$.

\demo Proof: Observe that $\g_\e=\{F_\e=0\} \cap co(k)$ where
$$
   F_\e(x_1,x_2) = x_1+x_2-\e(x_1^2+x_2^2)- 1+\e.
$$
Simple calculations show that for $0<\e< 1/2$ 
we have $F_\e(x_1,0)< 0$ for $0\le x_1<1$,
$F_\e(1,0)=F_\e(0,1)=0$, $F(x_1,1-x_1) > 0$ when $0< x_1< 1$, and
${\di \over \di x_2} F_\e(x_1,x_2) = 1-2\e x_2 >0$ for $0\le x_2\le 1$.
These properties imply that $\g_\e$ is a graph $y_1=h_\e(x_1)$
of a \ra\ function $h_\e$ over the segment $0\le x_1\le 1$, with
with the endpoints $(1,0)$ and $(0,1)$. Since
$\di F_\e/\di\e =1-(x_1^2+x_2^2)\ge 0$ on $co(k)$
we conclude that, as $\e$ decreases to $0$, the functions $h_\e$
increase to $h_0(x_1)=1-x_1$. This gives the last claim in lemma 4.3.

We will show that for sufficiently small $\e>0$ there exists a
bounded, simply connected region $D_\e \subset \{z_2=0\}$ with 
piecewise smooth boundary such that $\G_\e$ is the graph of 
a holomorphic function over $D_\e$. The equation for $\G_\e$ 
is equivalent to
$$ \eqalignno{ x_1+x_2-\e(x_1^2+x_2^2) +\e(y_1^2+y_2^2) &= 1-\e  \cr
   (1-2\e x_1)y_1 + (1-2\e x_2)y_2 &= 0.                          &(4.2) \cr}
$$
When $y_1=y_2=0$ we get the equation for $\g_\e$, and hence
$\G_\e \cap \R^2 = \g_\e$. On $co(k)$ we have $x_1+x_2\ge 0$
and $x_1^2+x_2^2 \le 1$, with equality only at the points $(1,0)$
and $(0,1)$. Rewriting the first equation in (4.2) in the form
$$
    (x_1+x_2) + \e(y_1^2+y_2^2) = 1-\e(1 -(x_1^2+x_2^2)) \le 1
$$
we see that (4.2) has no solutions for $|y|^2=y_1^2+y_2^2 > 1/\e$,
and it has no solutions on $\g_0 +i\R^2$ ($\g_0$ was defined
in lemma 4.3). Hence the boundary of $\G_\e$ is contained in
$K_\e$ and therefore $\G_\e \subset co(K_\e)$. 
From the second equation in (4.2) we get
$$
   y_2 = -y_1 {1-2\e x_1 \over 1-2\e x_2}                       \eqno(4.3)
$$
(again this requires $\e<1/2$ since $0\le x_1,x_2 \le 1$ on $\G_\e$).
Inserting this into the first equation (4.2) we get
$$
   G_\e(x_1,y_1,x_2) :=  x_1+x_2-\e(x_1^2+x_2^2) +
   \e y_1^2 \left( 1+ {(1-2\e x_1)^2 \over (1-2\e x_2)^2} \right) -1+\e =0.
                                                                \eqno(4.4)
$$
Consider first its restriction to $x_2=0$:
$$
    G_\e(x_1,y_1,0)=  x_1-\e x_1^2 +
   \e y_1^2 \left( 1+ (1-2\e x_1)^2 \right) -1+\e =0.
$$
Let $a_\e>0$ be the solution of the equation
$G(0,a_\e,0)= 2\e a_{\e}^2 -1+\e=0$. Calculations show that
$G_\e(0,y_1,0)<0$ for $|y_1|<a_\e$, $G_\e(1,y_1,0)\ge 0$
(with equality only at $y_1=0$), and
${\di G_\e \over\di x_1}(x_1,y_1,0)>0$ for $0\le x_1 \le 1$.
This shows that the set
$$
   \sigma_\e = \{x_1+iy_1 \colon 0 \le x_1\le 1,\ G_\e(x_1,y_1,0)=0 \}
$$
is a smooth \ra\ curve which can be written as a graph
$x_1= g_\e(y_1)$ over the interval $|y_1|\le a_\e$, and  the set
$$ \eqalign{
   D_\e &= \{x_1+iy_1 \in\C \colon 0< x_1< 1,\ G_\e(x_1,y_1,0)<0 \} \cr
        &= \{x_1+iy_1 \colon 0< x_1< g_\e(y_1),\ |y_1|< a_\e \} \cr}
$$
(with piecewise smooth boundary) is conformally equivalent to the disc.
A calculation shows that for $\e>0$ sufficiently small we have
${\di G_\e \over \di x_2}(x_1,y_1,x_2) >0$ on $0\le x_1\le 1$ and
$y_1^2\le 1/\e$, and $G_\e(x_1,y_1,1)>0$ for $x_1+iy_1\in D_\e$.
Since $G_\e(x_1,y_1,0)<0$ for $x_1+iy_1\in D_\e$, it follows that
(4.4) has a unique solution $x_2=\xi_\e(x_1,y_1) \in [0,1]$ for each
$z_1=x_1+iy_1\in \bar D_\e$ and it has no solutions for points in
$\{0\le x_1\le 1\}\bs \bar D_\e$. From (4.3) we also calculate $y_2$ and
thus obtain a unique analytic solution $z_2=f_\e(z_1)$
$(z_1\in\bar D_\e)$ of the system (4.2). This proves that $\G_\e$ is
an analytic disc with boundary in $K_\e$.
\endpr

We continue with the proof of proposition 4.2.
For each $y\in \R^2$ and $C \in \R$ we have
$$
    \{x\in\R^2 \colon \rho(x+iy)=C\} = \{(x_1,C)\colon x_1 \le C\}
     \cup \{(C,x_2)\colon x_2 \le C\}.
$$
For each point $z=x+iy\in\ctwo$ with
$C_0<\rho(z)<C_1$ we can choose a line segment
$l_z \subset R^2+iy$ passing through $z$ such that
$\rho>C_0$ on $l_z$ and the endpoints of $l_z$ 
belong to $\{\rho=C_1\}$. We can choose such $l_z$ 
depending smoothly on $z$ in the region $C_0<\rho(z)<C_1$.
The segment $l_z$ together with the two bounded segments in the 
level set $\rho=C_1$ (in $\R^2 +iy$) determines a closed triangle
$T_z \subset \R^2+iy$ which corresponds (after a rotation and dilation 
of coordinates) to the set $co(k)$ in the model case. 
Lemma 4.3, applied to a slightly larger triangle $\tilde T_z \supset T_z$ 
obtained by a small parallel translation of the segment 
$l_z$ so as to include the point $z$ in the interior of 
$\tilde T_z$, gives an analytic disc $\G_z \subset \ctwo$ 
passing through $z$ such that  $\rho>C_0$ on $\G_z$
and $\rho=C_1$ on $b\G_z$. We can parametrize $\G_z$ by a 
map $\l(z,\cdotp) \colon  \bar U \to \G_z$, holomorphic
in $U$ and depending continuously on $z\in K$.
\endpr

Combining proposition 4.2 an corollary 2.2 we obtain

\proclaim 4.4 Corollary: Let $\rho$ be the function (4.1).
Given a continuous map $g_0\colon  \bar U \to\ctwo$, \holo\ in $U$,
and constants $0<r<1$, $C_0,C_1\in \R$ that
$C_0 < \rho(g_0(\z)) < C_1$ for $r\le |\z|\le 1$, 
there is for each $\e>0$ a \hp\ map $g\colon  \bar U\to\ctwo$
satisfying
\item{(i)}   $|\rho(g(\z)) - C_1|< \e$\ \ $(|\z|=1)$,
\item{(ii)}  $\rho(g(\z))  > C_0$\ \ $(r\le |\z|\le 1)$, and
\item{(iii)} $|g(\z)- g_0(z)|<\e$ \ \ $(|\z|\le r)$.

\demo Proof of theorem 4.1:
Choose a sequence $\e_k>0$, $\sum_{k=1}^\infty \e_k <1$.
We begin by an arbitrary continuous map $g_1\colon  \bar U\to\ctwo$
that is holomorphic in $ U$ and a number $0<r_1<1$.
Choose numbers $M_0,M_1\in \R$ such that
$M_0 < \rho(g_0(\z)) < M_1$ for $r_1\le |\z|\le 1$.
Choose a number $M_2\ge M_1+1$ and apply corollary 4.4 to get
a polynomial map $g_2\colon \C\to\ctwo$ and a number $r_2$, $r_1<r_2<1$,
such that the following hold for $k=2$:
\item{(a$_k$)} $M_{k-1}< \rho(g_k(\z)) < M_k$\ \  $(r_k\le |\z|\le 1)$,
\item{(b$_k$)} $\rho(g_k(\z)) > M_{k-2}$\ \ $(r_{k-1}\le |\z|\le 1)$, and
\item{(c$_k$)} $|g_k(\z)-g_{k-1}(\z)| <\e_{k-1}$\ \ $(|\z|\le r_{k-1})$.

\ni This process can be continued inductively as follows. Suppose that we 
have already constructed $g_{k-1}$ for some $k\ge 2$. Choose 
$M_k\ge M_{k-1}+1$ and apply corollary 4.4 to get a map $g_k$ which 
satisfies (a$_k$) for $|\z|=1$ and it satisfies (b$_k$) and (c$_k$).
By continuity we can choose $r_k<1$ such that $1-r_k<(1-r_{k-1})/2$
and such that (a$_k$) holds for $r_k\le |\z|\le 1$.

By construction we have $\lim_{k\to\infty} r_k=1$, 
$\lim_{k\to\infty} M_k =+\infty$,
and $g=\lim_{k\to\infty} g_k$ exists uniformly on compacts in $ U$ by (c$_k$).
It remains to show that $\rho(g(\z)) \to\infty$ as $|\z|\to 1$. Fix $k\ge 2$
and consider points in $A_k=\{ \z\colon r_{k-1}\le |\z|\le r_k\}$.
For $l\ge k$ we have $|g_{l+1}(\z)-g_l(\z)| < \e_l$, so
$|g(\z)-g_k(\z)| \le \sum_{l=k}^\infty |g_{l+1}(\z)-g_l(\z)| <
                     \sum_{l=k}^\infty \e_l <1.
$
From this and (b$_k$) we get
$\rho(g(\z)) > \rho(g_k(\z)) -1 > M_{k-2}-1$ for $\z\in A_k$.
Since $M_{k-2}\to\infty$ as $k\to\infty$, the result follows.
\endpr

\demo Remark: One can give an alternative proof of theorem 1.4
as follows. One can construct a family of \holo\ maps 
$F_p\colon \C\to\C^2$, depending continuouly on 
$p\in (\C^*)^2$, such that (i) $F_p(0)=p$,
(ii) $|F_p(\z)|\ge |p|-\e_p$ for all $\z\in \C$
(where $\e_p>0$ can be made independent of $p$ in any
compact set $K\ss(\C^*)^2$, (iii) $F_p(\C)$ misses $zw=0$, 
and (iv) $\lim_{|\z|\to\infty} |F_p(\z)|=+\infty$.
The discs $\z\to F(\z)$, $|\z|\le R$ with $R$ large enough, 
can be taken as building blocks to construct proper
\holo\ discs $U\to\C^2$ whose image avoids both coordinate
axes (compare with proposition 4.2). Similar method can
used to construct proper \holo\ discs in $\C^2$ avoiding 
the curve $zw=1$.

%
%
%

\beginsection 5. Boundary behavior of proper holomorphic discs.

In this section we prove theorem 1.5. We begin by recalling 
some classical results on boundary behavior of meromorphic 
functions on $U=\{|\z|<1\}$ (see e.g.\ [CL] and [Pri]). 
Let $\bar \C=\C\cup\{0\}$ denote the Riemann sphere. 
For $a\in \C$ and $r>0$ set $D(a;r)=\{\z\in\C\colon |\z-a|<r\}$.
In what follows let $f$ be a meromorphic function on $U$.
We denote by $C(f,e^{i\theta})$ its {\bf unrestricted
cluster set} at $e^{i\theta}\in T$:
$$
    C(f,e^{i\theta})=
    \bigcap_{r>0} \bar{f\bigl( U\cap D(e^{i\theta};r)\bigr)}.
$$
Equivalently, $a\in\bar\C$ belongs to $C(f,e^{i\theta})$ \iff\
there exists a sequence $\z_j\in  U$ such that
$\lim_{j\to\infty} \z_j = e^{i\t}$ and $\lim_{j\to\infty} f(\z_j)=a$.
If $D\subset  U$ is a subset with $\eit\in\bar D$, we denote by
$C_D(f,\eit)$ the {\bf restricted cluster set} of $f$ at $\eit$,
defined as the set of limits of $f$ along sequences $\z_j\in D$
with $\lim_{j\to\infty} \z_j = \eit$.

A point $e^{i\theta}$ for which $C(f,e^{i\theta})= \bar \C$ is called
a {\bf Weierstrass point} of $f$, and the set of all such points is
the {\bf Weierstrass set} $W(f)$ [CL, p.\ 149].

For each $\e^{i\t} \in T$ and $0<\a<1$ we set
$$
    \G_\a(e^{i\t}) = \{\z\in U\colon
    |\Im (1-\z e^{-i\t})| < \a |\z - e^{i\t}| \}.           
$$
This is an angle in $U$ with vertex at $\eit$ and opening
$2\arcsin \a$, bisected by the radius that terminates at $\eit$.
If the limit
$$
    f^*(\eit)= 
	\lim_{\G_\a(\eit) \owns  \z\to\eit} f(\z) \in \bar \C  \eqno(5.1)
$$
exists and is independent of $\a$, it is called the 
{\it nontangential limit of $f$ at $\eit$} and
$\eit$ is called a {\bf Fatou point} of $f$. The set of 
all Fatou points the {\bf Fatou set} $F(f)$ [CL, p.\ 21].

A point $\eit \in T$ is called a {\bf Plessner point}
of $f$ if for every angle $\G$ with vertex $\eit$ the
partial cluster set $C_\G(f,\eit)$ equals $\bar \C$ 
(i.e., it is total). The set of all Plessner points is 
the {\bf Plessner set} $I(f)$ [CL, p.\ 147].
Clearly $I(f)\subset W(f)$. 

The {\bf Nevanlinna characteristic} of a \holo\ 
function $f$ on $U$ is defined by 
$$ T(r,f) = \int_0^{2\pi} \log^+ |f(r\eit)|\, {dt\over 2\pi} 
   \qquad (0\le r<1); 
$$
for meromorphic functions see [CL, p.\ 39] or [Nev].

The {\bf range} of $f$, denoted $R(f)$, is the set 
of all $\a\in \bar \C$ for which there exists a 
sequence $\z_j	\in U$ with $\lim_{j\to \infty}\z_j =1$ 
such that $f(\z_j)=\a$ for all $j\in \N$.
By restricting the attention only to sequences 
$z_j\in U$ with $\lim_{j\to\infty} \z_j =\eit$ we 
get the range of $f$ at $\eit$, denoted $R(f,\eit)$.

The notion of {\bf logarithmic capacity} of a Borel set
$E\subset \C$ can be found in [CL, p.\ 9]. Such a set is 
of capacity zero \iff\ it is polar, i.e., it is contained 
in the $-\infty$ level set of a non-constant subharmonic 
function on $\C$ [Lan, Tsu]. The following summarizes some of 
the known results which we shall need in the proof of theorem 1.5.

\medskip\ni\bf 5.1 Theorem: \sl 
Let $f$ be a meromorphic function on the disc $U$.
\item{(a)} If $\eit\in T$ is not a Weierstrass point of $f$
then there is an open arc $\g\subset T$ containing $\eit$
such that almost every point in $\g$ is a Fatou point of $f$.
\item{(b)} If $f$ has bounded Nevanlinna characteristic on $U$
then almost every point of $T$ is a Fatou point of $f$.
\item{(c)} Almost every point in $T$ belongs to $F(f)\cup I(f)$.
\item{(d)} If $f(U)\subset \bar\C\bs E$ for some set
$E$ of positive capacity then $f$ has bounded Nevanlinna 
characteristic. 
\item{(e)} If $R(f,\eit)$ omits a set $E\subset \bar \C$ 
of positive capacity then there is an open arc 
$\gamma\subset T$ containing $\eit$ such that
almost every point of $\gamma$ is a Fatou point of $f$.
\rm

\demo Proof:
(a) If $\eit$ is not a Weierstrass point of $f$,
there is a disc $D(\eit;r)$ such that
$f(D(\eit;r)\cap U)$ omits a disc $D(a;\d) \subset\C$.
The function $g(\z)=1/(f(\z)-a)$ is then bounded
holomorphic in $D(\eit;r)\cap  U$ and hence 
by Fatou's theorem it has nontangential limit 
$g^*(e^{it})$ at almost every point 
$e^{it}\in \g= T\cap D(\eit;r)$ [CL, p.\ 21]. 
The same is then true for $f$ and hence almost every 
point of $\g$ belongs to the Fatou set $F(f)$. 
(See also [CL, Theorem 8.4].)
Part (b) follows by combining Fatou's theorem
with a theorem of R.\ Nevanlinna to the effect that
a meromorphic function with bounded Nevanlinna
characteristic on $U$ is the quotient of two bounded
\holo\ functions [CL, p.\ 41].
Part (c) is a classical theorem due to Plessner 
([Ple], [Pri, p.\ 217] or [CL, p.\ 147]).

Part (d) is due to Frostman [Fro]; the following 
simple proof was shown to us by D.\ Marshall.
After a fractional linear transformation we may 
assume that $\infty\in E \subset \{|z|>1\}$. 
Let $g_0(z)$ be the Green's function for $\C\bs E$ with 
a logarithmic pole at $0$ (so $g_0(z)\to 0$ as $z\to E$).
Then $\log^+{1\over |z|} \le g_0(z)$. The function 
$u(z):= g_0(z) +  \log|z|$ is harmonic on $\C\bs E$ 
since both summands are harmonic on $\C\bs (E\cup\{0\})$ 
and the pole at $0$ cancels off.
If $f\colon U\to \C\bs E$ is a \holo\ function then
$u\circ f$ is harmonic on $U$ and we have
$$ \eqalign{ \int_0^{2\pi} \log^+  |f(r\eit)|\, {d\theta\over 2\pi}
   &= \int_0^{2\pi} \left( 	\log^+ {1\over |f(r\eit)|}
     +  \log{|f(r\eit)|} \right) \,{d\theta\over 2\pi} \cr
	&\le  \int_0^{2\pi} \left( g_0(f(r\eit)) +
     \log{|f(r\eit)|} \right) \, {d\theta\over 2\pi} \cr
    &= \int_0^{2\pi} u(f(r\eit)) \, {d\theta\over 2\pi}
	 = u(f(0)). \cr}
$$
Thus $\int_0^{2\pi} \log^+ |f(r\eit)|\, {dt\over 2\pi} \le u(f(0))$
for $r\in (0,1)$ which proves (d).

For (e) observe that $R(f,\eit)=\bigcap_{n\in\N} f(D_n)$, 
where $D_n=D(\eit;1/n)\cap U$. The sets $f(D_n)$ are decreasing
with $n$. If $R(f,\eit)$ omits a set $E$ of positive capacity 
then $f(D_n)$ omits a set $E'$ of positive capacity 
for some sufficiently large $n\in \N$. We may assume that 
$\infty\in E'$. Observe that $D_n$ is conformally equivalent
to the disc. From (d) and (b) applied to the holomorphic 
function $f\colon D_n\to \C\bs E'$ it follows that almost 
every point of the arc $\gamma= D_n \cap T$ is a 
Fatou point of $f$.
\endpr

We shall frequently use the following uniqueness
theorem due to Plessner ([Ple], [CL, p.\ 146])
and to Lusin and Priwalow [Pri, p.\ 212].

\proclaim 5.2 Theorem:
If a meromorphic function $f$ on $U$ has an angular limit
equal to zero at each point in a set $E\subset T$
of positive Lebesgue measure then $f$ is the zero function.

\demo Remark: In theorem 5.2 {\it we cannot replace
angular limits with radial limits}, see examples due to 
Lusin and Priwalow in [Pri], sec.\ IV.5. Here we use the 
term `angular limit' rather than `nontangential limit'
since the latter usually means the existence of the limit 
within every angle with the given vertex.

\demo Proof of theorem 1.5: Let
$(f_1,f_2)\colon  U\to\ctwo$ be a \phm\ and
let $g$ be any of the functions as in theorem 1.5.
It suffices to show that the Fatou set $F(g)$ has measure 
zero. From theorem  5.1 (a) it will then follow that
$W(g)=T$, theorem 5.1 (c) will imply that
the Plessner set $I(g)$ has full measure in $T$,
and theorem 5.1 (e) will imply that the complement
of the range $R(g,\eit)$ in $\C$ has capacity zero
for each $\eit\in T$. Since sets of capacity zero in 
$\C$ coincide with polar sets ([Tsu], [Lan]), 
theorem 1.5 (d) follows. 

To prove that $F(g)$ has measure zero we 
consider separately each case.

\demo Case (i): 
Suppose that $f_1$ has an angular limit $f_1^*(\eit) \in\bar\C$
(5.1) at all points $\eit$ in a set $A\subset T$.
Then $A$ is Lebesgue measurable and can be written as
$A=A_1\cup A_2$, where $A_1$ is the set of all $\eit \in A$ such that
$f_1^*(\eit)\in \C$ and $A_2$ is the set of all $\eit \in A$ with
$f_1^*(\eit)=\infty$. Then $1/f_1$ has angular limits zero at each
point of $A_2$. If $A_2$ is of positive measure,
theorem 5.2 implies that $1/f_1$ is identically zero
in $U$, a contradiction.  Thus $A_2$ has measure zero.
Consider now $A_1$. Since $(f_1,f_2)\colon U\to\ctwo$ is proper,
$\max\{|f_1(\z)|,|f_2(\z)|\}$ tends to $+\infty$ as $|\z|\to 1$.
Since $f_1$ has a finite angular limit at each $\eit\in A_1$,
$|f_2|$ has an angular limit $\infty$ at each point of $A_1$.
If $A_1$ is of positive measure, Plessner's theorem, applied to
$1/f_2$, gives a contradiction as before. This shows that $A_1$
is of measure zero as well, and therefore the Fatou set of $f_1$
is of measure zero. The same applies to $f_2$.

\demo Case (ii):
Suppose that $g=f_1/f_2$ has an angular limit $g^*(\eit) \in\bar \C$
(5.1) within an angle $\G_\t$ at each point $\eit$ in a set
$A\subset T$. As in part (i) we write $A=A_1\cup A_2$,
where $g^*$ is finite on $A_1$ and equals $\infty$ on $A_2$.
Theorem 5.2 shows as above that $A_2$ must be of measure zero
for otherwise $g$ would be constant. If $A_1$ is of positive measure,
there is a set $A_0\subset A_1$ of positive measure and a
number $0<M<\infty$  such that $|g^*(\eit)| <M$ for each $\eit\in A_0$.
Hence there is a disc $U_\t$ centered at $\eit$ such that
$|f_1(\z)/f_2(\z)| \le M$ for $\z \in \G_\t \cap U_\t$.
Hence $|f_1(\z)| \le M |f_2(\z)|$ and therefore
$$
    \max\{|f_1(\z)|,|f_2(\z)|\} \le
    \max\{M|f_2(\z)|,|f_2(\z)|\}
                                \quad (\z \in \G_\t \cap U_\t).
$$
Since this maximum tends to $+\infty$ as $\z\to\eit$, it follows that
$|f_2(\z)| \to\infty$ as $\z\to\eit$ within $\G_\t$. Thus $1/f_2$ has angular
limits zero at each point of $A_0$, a contradiction to theorem 5.2.
This proves that $A_1$ must be of measure zero as well.

\demo Case (iii):
This follows from case (i) by observing that
for each nonconstant \holo\ polynomial $P$ on $\ctwo$ \te\
another \holo\ polynomial $Q$ such that $(P,Q)\colon \ctwo\to\ctwo$
is a proper map, and hence $P(f_1,f_2) \colon  U \to \ctwo$
is the first component of a proper map $U \to\ctwo$.
In fact, we have

\proclaim 5.3 Lemma: Let $P$ and $Q$ be nonconstant \holo\ polynomials
on $\ctwo$ whose leading order homogeneous parts $P'$ resp.\ $Q'$ have
no common zero on $\ctwo\bs\{0\}$. Then $(P,Q)\colon \ctwo\to\ctwo$
is a proper map.

We leave out the simple proof. Observe that the zero set of $P'$
is a finite union of complex lines, so it suffices to choose $Q$
to be a linear function which does not vanish on $P'=0$ except at
the origin; the pair $(P,Q)$ then provides a proper self-map of $\ctwo$.

\demo Case (iv): 
Apply (i) and lemma 5.3 to the map
$\bigl(P(f_1,f_2),Q(f_1,f_2)\bigr) \colon  U \to\ctwo$.

\medskip\ni\bf Appendix: 
Crossing a critical level by analytic discs. \rm

Let $X$ be a Stein manifold of dimension at least two
and $\rho\colon X\to\R$ a \spsh\ Morse exhaustion function.
Let $p\in X$ be a critical point of $\rho$. Choose 
constants $c_0,c_1$ such that $c_0 <\rho(p)<c_1$ 
and $p$ is the only critical point of $\rho$ 
in $A(c_0,c_1)=\{x\in X\colon c_0\le \rho(x)\le c_1\}$.
Suppose that $f_0\colon \bar U\to X$ is a \holo\ map
such that $c_0 < \rho(f_0(\eit)) <\rho(p)$ 
for each $\eit\in T$. In [Glo]  the second author
showed how to construct a smooth map $f_1\colon \bar U \to X$
which is close to being holomorphic on $U$ such 
that $\rho(p)< \rho(f_1(\eit)) < c_1$ for $\eit\in T$
and such that $f_1$ approximates $f_0$ on a smaller
disc $|\z|\le r<1$. The map $f_1$ is obtained by adding
to $f_0$ a small non-holomorphic contribution which can
be controlled by the data. Once the boundary curve $f_1(T)$
passes the critical level of $\rho$ at $p$ we can 
use the procedure described in sect.\ 2 above (or in [FG])
to continue pushing it higher towards the next critical 
level of $\rho$. It was shown in [Glo] that the non-holomorphic 
contribution made at the initial step may be cancelled off during 
a later stage of the construction, once the boundary of the 
disc is sufficiently far above the critical level at $p$. 
The reason is that the modification process is a linear one,
and we obtain the final solution as the sum of a 
convergent series. (Here it is convenient to embed
$X$ into a Euclidean space $\C^N$.) 

Here we wish to point out that the transition from 
$f_0$ to $f_1$ as above can also be accomplished by applying 
to $f_0$ the gradient flow $\theta_t$ of $\rho$ (in the direction
of increasing $\rho$). Unless a point $x\in A(c_0,c_1)$ belongs 
to the stable manifold $W^s(p)$ of $p$ (see e.g.\ Shub [Sh]), 
we have $\rho(\theta_t(x)) > \rho(p)$ for sufficiently large
$t>0$. Thus, if $f_0(T)\cap W^s(p)=\emptyset$, we
can choose a smooth positive function $a$ on $\bar U$
such that the map $f_1(\zeta) =\theta_{a(\z)}(f_0(\z))$
($|\z|\le 1$) satisfies $\rho(f_1(\eit))>\rho(p)$
for $\eit\in T$. 

If the number $c_0$ is sufficiently close to $\rho(p)$ 
as we may assume to be the case, the set $W^s(p) \cap A(c_0,c_1)$ 
is a closed real submanifold of $A(c_0,c_1)$ 
whose dimension equals the index $i(p)$ 
(the number of negative eigenvalues of the Hessian) 
of $\rho$ at $p$. Since $\rho$ is \spsh, we have 
$i(p) \le \dim_{\C} X$ (see [AF]) and therefore 
$\dim f_0(T)+\dim W^s(p)\le 1+ \dim_{\C} X <\dim_{\R}X$.
By transversality a generic small \holo\ perturbation of 
$f_0$ satisfy the required condition 
$f_0(T)\cap W^s(p)=\emptyset$ which makes it possible
to obtain $f_1$ as above. The rest of the procedure 
remains as in [Glo].

%
%
%
%
\bigskip
\ni\bf References. \rm \medskip

\item{[AF]} A.\ Andreotti, T.\ Frankel:
The Lefschetz theorem on hyperplane sections. 
Ann.\ of Math.\ (2) {\bf 69}, 713--717 (1959). 

\item{[CL]} E.\ F.\ Collingwood, A.\ J.\ Lohwater:
The theory of cluster sets.
Cambridge Tracts in Mathematics and Mathematical Physics {\bf 56},
Cambridge University Press, Cambridge, 1966.

\item{[FG]} F.\ Forstneri\v c, J.\ Globevnik:
Discs in pseudoconvex domains.
Comment.\ Math.\ Helv.\ {\bf 67} 129--145 (1992).

\item{[Fro]} O.\ Frostman: Potentiel d'equilibre et capacit\'e
des ensembles avec quelques applications \`a la th\'eorie
des fonctions.
Medd.\ Lunds Mat.\ Sem.\ {\bf 3}, 1--118 (1953).

\item{[Glo]} J.\ Globevnik: Discs in Stein manifolds.
Indiana Univ.\ Math.\ J.\ {\bf 49} (2000), 553--574.

\item{[H\"or]} L.\ H\"ormander:
An Introduction to Complex Analysis in Several Variables.
Third edition. North-Holland, Amsterdam-New York, 1990.

\item{[Lan]}  N.\ S.\ Landkof:
Foundations of Modern Potential Theory.
Translated from the Russian by A.\ P.\ Doohovskoy.
Die Grundlehren der mathematischen Wissenschaften,
Band 180. Springer, New York--Heidelberg, 1972.

\item{[Nev]} R.\ Nevanlinna:
Eindeutige analytische Funktionen, 2te Aufl.,
Springer, Berlin--G\"ot\-tingen--Heildeberg 1953.

\item{[Ple]} A.\ I.\ Plessner:
\"Uber das Verhalten analytischer Funktionen
am Rande ihres Definitionsbereiches.
J.\ reine angew.\ Math.\ {\bf 158}, 219--227 (1927).

\item{[Pri]} I.\ I.\ Priwalow:
Randeigenschaften analytischer Funktionen, 2.ed.
Hoch\-schulb\"ucher f\"ur Mathematik, Bd.\ 25.
VEB Deutscher Verlag der Wissen\-schaften, Berlin, 1956.

\item{[Sh]}  M.\ Shub: Global stability of dynamical systems.
Translated from the French by Joseph Christy. 
Springer, New York-Berlin, 1987.

\item{[SW]} N.\ Sibony, P.\ M.\ Wong:
Some results on global analytic sets.
S\'eminaire Pierre Lelong--Henri Skoda
(Analyse). Ann\'ees 1978/79 (French), pp.\ 221--237,
Lecture Notes in Math., 822, Springer, Berlin, 1980.

\item{[Tsu]} M.\ Tsuji:
Potential Theory in Modern Function Theory.
Reprinting of the 1959 original. Chelsea
Publishing Co., New York, 1975.

%
%
%
%
\bigskip\medskip
\itemitem{\it Address:} Institute of Mathematics, Physics and Mechanics,
University of Ljub\-ljana, Jadranska 19, 1000 Ljubljana, Slovenia

\bye